\documentclass[smallextended]{svjour3}       
\smartqed  
\usepackage{graphicx}
\usepackage[utf8]{inputenc} 
\usepackage[T1]{fontenc}    
\usepackage{hyperref}       
\usepackage{url}            
\usepackage{microtype}
\usepackage{array}
\usepackage{color} 
\usepackage{tabularx}
\usepackage{amsmath,mathdots}
\usepackage{amssymb}
\usepackage{amsfonts}
\usepackage{txfonts}
\usepackage{xcolor}
\usepackage{bm}
\usepackage{soul}
\usepackage{float}
\usepackage{cite}
\usepackage{booktabs, multirow}

\hypersetup{
	colorlinks=true,                          
	linkcolor=blue, 
	citecolor=blue, 
	urlcolor=black  
} 
\usepackage[colorinlistoftodos]{todonotes}

\hypersetup{
	pdftitle={},
	pdfsubject={},
	pdfauthor={Alessia Lucca, Michael Dumbser},
}

\usepackage{xcolor}
\definecolor{darkgreen}{rgb}{0.1,0.6,0.1}

\newcommand{\halb}{\frac{1}{2}}

\newcommand{\A}{\mathbf{A}} 
\newcommand{\B}{\mathbf{B}} 
\newcommand{\E}{\mathbf{E}} 
\newcommand{\D}{\mathbf{D}} 

\renewcommand{\H}{\mathbf{H}}

\newcommand{\p}{\mathbf{p}} 
\newcommand{\q}{\mathbf{q}} 
\newcommand{\f}{\mathbf{f}} 
\newcommand{\x}{\mathbf{x}} 
\newcommand{\pd}{\partial} 
 
\newcommand{\err}{\mathit{r}} 
\newcommand{\normal}{\mathbf{n}} 

\newcommand{\tL}{\tilde{L}}
\newcommand{\tp}{\tilde{\mathbf{p}}} 
\newcommand{\tH}{\tilde{\mathbf{H}}} 
 
\newcommand{\tE}{\tilde{\mathbf{E}}} 
 
\newcommand{\vv}{\mathbf{v}} 
\newcommand{\uu}{\mathbf{u}} 
 
\newcommand{\tuu}{\tilde{\mathbf{u}}} 
\newcommand{\ts}{\tilde{s}} 
\newcommand{\en}{\mathcal{E}}

\allowdisplaybreaks

\begin{document}
	
	\title{Structure-preserving schemes for nonlinear symmetric hyperbolic and thermodynamically compatible systems of partial differential equations}
	
	\titlerunning{Structure-preserving schemes for nonlinear SHTC systems}        
	
	\author{Alessia Lucca \and 
		Michael Dumbser
	}
	
	\institute{Alessia Lucca  \at 
		Department of Civil, Environmental and Mechanical Engineering, University of Trento, Via Mesiano 77, 38123 Trento, Italy.\\
		\email{alessia.lucca@unitn.it}
		\and
		Michael Dumbser \at
		Department of Civil, Environmental and Mechanical Engineering, University of Trento, Via Mesiano 77, 38123 Trento, Italy.\\
		\email{michael.dumbser@unitn.it}
	}
	
	\date{Received: December 2025}
	
	\maketitle


\begin{abstract}
This paper aims at developing exactly energy-conservative and structure-preserving finite volume schemes for the discretisation of first-order symmetric-hyperbolic and thermodynamically compatible (SHTC) systems of partial differential equations in continuum physics. Due to their thermodynamic compatibility the class of SHTC systems satisfies an additional conservation law for the total energy and many PDE in this class of equations also satisfy stationary differential constraints (involutions). First, we propose a simple semi-discrete cell-centered HTC finite volume scheme that employs collocated grids and that is compatible with the total energy conservation law, but which does not satisfy the involutions. Second, we develop a fully discrete semi-implicit finite volume scheme that conserves total energy and which can be proven to satisfy also the involution constraints exactly at the discrete level. This method is a vertex-based staggered semi-implicit scheme that preserves the basic vector calculus identities $\nabla \cdot \nabla \times \A = 0$ and $\nabla \times \nabla \phi = 0$ for any vector and scalar field, respectively, exactly at the discrete level and which is also exactly totally energy conservative. The main key ingredient of the proposed implicit scheme is the fact that it uses a discrete version of the \textit{symmetric-hyperbolic Godunov-form} of the governing PDE system. This leads naturally to sequences of symmetric and positive definite linear algebraic systems to be solved inside an iterative fixed-point method used in each time step. 
We apply our new schemes to three different SHTC systems. In particular, we consider i) the equations of nonlinear acoustics, ii) the nonlinear Maxwell equations in the absence of charges, and iii) a nonlinear version of the Maxwell-GLM system. We also show some numerical results to provide evidence of the stated properties of the proposed schemes. 
\end{abstract}


\section{Introduction}

Structure-preserving (SP) numerical methods are designed to preserve, even at a discrete level, some intrinsic physical properties, geometric characteristics or mathematical structures of the continuous model during the simulation process. Besides the conservation of typical quantities such as mass, momentum and energy, a wide field in research is also devoted to the development of numerical schemes which preserve inherent mathematical features of the system being modeled, such as stationary differential constraints, so-called involutions.  
Over time several hyperbolic PDE systems proposed for the description of dynamic processes in continuum physics are endowed with involution constraints. Involution constraints in general are stationary differential equations that are satisfied by the governing PDE system at all times if they are satisfied by the initial data. The most prominent example is the system of the Maxwell equations, where the magnetic field must remain divergence-free as well as the divergence of the electric field vanishes in absence of charges if they were initially. Another famous involution is the curl-free condition of the velocity field in the system of the acoustic wave equations.

Several numerical frameworks have been developed to guarantee the divergence-free property exactly at the discrete level. Some make use of so-called mimetic finite differences \cite{Yee66,HymanShashkov1997,Margolin2000,Lipnikov2014} and define the evolution quantities on appropriately staggered grids, some have developed compatible finite volume schemes \cite{DeVore,GardinerStone,BalsaraCED,BalsaraKaeppeli,HazraBalsara}, while others compatible finite elements \cite{Nedelec1,Nedelec2,Hiptmair,Arnold2006,Monk,Alonso2015,Zampa1,Zampa2}. 
Very recently a new framework of discontinuous Galerkin finite element schemes on unstructured meshes has been proposed which, in suitable combination with continuous Lagrange finite elements, achieve an exact discrete compatibility with the curl and divergence involutions \cite{abgrall2025simple}. The aforementioned method is a direct extension to higher order in space of the vertex-based compatible discretizations introduced in \cite{Maire2007,Maire2020,Barsukow2024,Sidilkover2025}. Alternative involution-preserving DG schemes have been recently introduced in \cite{Perrier1,Perrier2,Perrier3}.  
On the other hand, more rare are curl-preserving numerical schemes for PDE with curl involutions, see e.g. \cite{Torrilhon2004,JeltschTorrilhon2006,BalsaraCurlFree,SIGPR,SIST,Dhaouadi2023NSK}.
Rather than enforcing divergence constraints on the electric and the magnetic field exactly, a new approach has been proposed in Munz et al. \cite{MunzCleaning}, called generalized Lagrangian multiplier (GLM) divergence cleaning. Here, the core concept is to solve an augmented evolution system, where artificial scalar cleaning variables are incorporated and linked to the original equations. This ensures that divergence errors in the magnetic and electric field are unable to accumulate, but instead propagate away through acoustic-type waves. A compareble strategy has been recently proposed also for the numerical handling of curl inaccuracies in curl-free vector fields, see approach has recently also been forwarded for the numerical treatment of curl errors in curl-free vector fields, see \cite{dumbser2020glm,SHTCSurfaceTension,HyperbolicDispersion,Dhaouadi2022,TwoFluidDG,FirasLauraCurlGLM}.

All of the aforementioned mathematical models fall into the larger class of symmetric hyperbolic and thermodynamically compatible (SHTC) systems, introduced by Godunov and and Romenski in \cite{God1961,Rom1998} and subsequently extended by them and their collaborators in \cite{GodunovRomenski72,GodRom2003,Rom1998,RomenskiTwoPhase2010,Godunov2012,GRGPR}.
The SHTC equations are derived from the minimization of the action formed by a Lagrangian. The involution constraints appear intrinsically in the Lagrangian form of SHTC equations as integrability conditions. A essential feature of SHTC systems is that they satisfy an additional conservation law for the total energy density, which allows the systems to be transformed into a symmetric form by introducing a generating potential and thermodynamic dual variables, the so-called main field, see \cite{Ruggeri81}. Furthermore, if the generating potential is a convex function, the SHTC system is symmetric hyperbolic in the sense of Friedrichs \cite{FriedrichsSymm}. Therefore, in the SHTC system the total energy conservation law is not a primary evolution equation. It is rather obtained as a linear combination of all the other equations, due to the privileged role that the total energy plays in the underlying variational principle from which all SHTC systems can be derived.

The proposed work focuses on three types of nonlinear SHTC system with involutions: the nonlinear acoustic equations, the nonlinear Maxwell system, and a nonlinear version of the Maxwell-GLM system. In these systems, a general relation is considered between the main variables and their associated thermodynamic dual variables, which may be nonlinear. For these systems two compatible discretizations are proposed. The first scheme is a semi-discrete finite volume method on collocated grids that preserves the total energy exactly. The second scheme is a vertex-based staggered semi-implicit scheme that preserves the basic vector calculus identities $\nabla \cdot \nabla \times \A = 0$ and $\nabla \times \nabla \phi = 0$ exactly on the discrete level and which is also totally energy conservative at the fully discrete level.

The rest of this paper is organized as follows. In Section \ref{sec.model}, we present the continuous models studied in this work and we will show their symmetric hyperbolicity. In particular, we will consider the nonlinear acoustic equations, the nonlinear Maxwell equations in vacuum and the augmented Maxwell-GLM system extended to the case of a general nonlinear Lagrangian. In Section \ref{sec.scheme}, we construct two compatible discretisations. The first one is a cell-centered semi-discrete thermodynamically compatible finite volume scheme that satisfies the extra conservation law for total energy at the discrete level. The second one is a staggered fully discrete, semi-implicit, finite volume scheme that is totally energy conservative and achieves discrete compatibility with curl and divergence involutions with the aid of suitable mesh staggering and mimetic finite difference operators. Some numerical results are presented in Section \ref{sec.results}.
The conclusions and an outlook to future work are given in Section \ref{sec.Conclusions}. 

\section{The continuous models}
\label{sec.model}
All overdetermined systems of nonlinear hyperbolic, thermodynamically compatible partial differential equations considered in this paper and detailed below, can be cast into the general form of a nonlinear hyperbolic system of conservation laws 
\begin{equation}
	\frac{\partial \q}{\partial t} + \frac{\partial \f_k(\q)}{\partial x_k}=0
	\label{eqn.pde}
\end{equation}
and they are all endowed with an extra conservation law for the total energy density $\en$ of the form
\begin{equation}
	\frac{\partial \en}{\partial t} + \frac{\partial F_k}{\partial x_k}=0.
	\label{eq:conE}
\end{equation}
Here, $ \q $ and $ \f_k(\q) $ are the state vector and the flux tensor, respectively, whereas in the additional conservation law \eqref{eq:conE}, the total energy flux is denoted by $ F_k=F_k(\q) $.
According to the formalism introduced by Godunov and Romenski in \cite{God1961,GodunovRomenski72,Rom1998} and Boillat in \cite{Boillat74}, if we denote the state vector of main evolution variables by $\q $, then there is a vector of associated thermodynamic \textit{dual variables}, also called the \textit{main field} \cite{Ruggeri81} defined as
\begin{equation}
	\p = \frac{\partial \mathcal{E}}{\partial \q} = \mathcal{E}_{\q}.
	\label{eqn.pdef} 
\end{equation}
The total energy density $\mathcal{E}$ is then linked to the so-called \textit{generating potential} $L$ via the Legendre transform
\begin{equation}
	L = \mathbf{p} \cdot \mathbf{q} - \mathcal{E} = \mathcal{E}_{\q} \cdot \q - \mathcal{E}. 
	\label{eqn.legendre} 
\end{equation}
It is easy to see that \eqref{eqn.pdef} and \eqref{eqn.legendre} imply
\begin{equation}
	\q = \frac{\partial L}{\partial \p} = L_{\p} \quad \text{and} \quad L_{\p\p} = \frac{\partial L_\p}{\partial p}= \frac{\partial \q}{\partial \en_\q} = \left(\en_{\q\q}\right)^{-1}.   	
\end{equation}
Thanks to the compatibility of \eqref{eqn.pde} with the extra conservation law for the total energy \eqref{eq:conE},
the following identity must hold:
\begin{equation}
	\p \cdot \frac{\partial \mathbf{f}_k(\q)}{\partial x_k} = \frac{\partial F_k}{\partial x_k}.   
	\label{eqn.flux.comp} 
\end{equation}
Throughout this paper we adopt the Einstein index notation assuming summation over two repeated indices.

\paragraph{The system of nonlinear acoustics.}

The governing equations for nonlinear acoustics read
\begin{align}
	\frac{\partial \vv}{\partial t} + \nabla s  &=0, \label{eq:Acoustic.v}\\
	\frac{\partial \rho}{\partial t} + \nabla \cdot \uu &=0, \label{eq:Acoustic.rho}
\end{align}
and satisfy the involution $ \nabla\times\vv=0 $. 
We identify the field $ \vv $ and the scalar $ \rho $ as state variable, whereas 
$ \uu = \partial \en/\partial \vv $ and $ s = \partial \en /\partial \rho $ represent their corresponding dual variables. 
Multiplying \eqref{eq:Acoustic.v}--\eqref{eq:Acoustic.rho} with the dual variables $ \p $ and summing up
leads the following total energy conservation law
\begin{equation}
	\frac{\pd \en}{\pd t} + \nabla \cdot \left(\uu s\right) =0.
\end{equation}
As in the linear case, the nonlinear acoustic system can also be derived from variational principles. This can be achieved by formally introducing a scalar potential $ Z $ as follows
\begin{equation}
	\vv = \nabla Z, \qquad s = - \partial t Z
	\label{eqn.dynamic1} 
\end{equation}
and by reparametrizing the Lagrangian according to 
\begin{equation}
	\mathcal{L}(\vv,s) =  \Lambda(Z, \pd_t Z, \nabla Z ).
\end{equation}
The first acoustic equation \eqref{eq:Acoustic.v} is then obtained as direct consequence of the definitions \eqref{eqn.dynamic1} by differentiating the first relation with respect to time and using the second identity. Introducing the relations
\begin{equation}
	\frac{\pd \Lambda}{\pd \partial_t Z} =-\frac{\pd \mathcal{L}}{\pd s}, \qquad \frac{\pd \Lambda}{\pd \partial Z_k} =\frac{\pd \mathcal{L}}{\pd v_k}
\end{equation}
and assuming $\pd \Lambda/\pd Z = 0$, the variation with respect to $Z$ gives the following Euler-Lagrange equations
\begin{equation}
	\frac{\pd \mathcal{L}_{s}}{\pd t}  -  \pd_k \mathcal{L}_{v_k} = 0.
	\label{eq:EL.s} 
\end{equation}
We then introduce a partial Legendre transform to $\mathcal{L}(\uu, s)$ and the new state variables 
\begin{equation}
	\rho := \mathcal{L}_{s}, \qquad
	\en(\vv,\rho) := s \mathcal{L}_{s} - \mathcal{L},
	\label{eq:rel1}
\end{equation}
with the following relations 
\begin{equation}
	\uu = \en_{\vv} = -\mathcal{L}_{\vv}, 
	\qquad
	s= \en_{\rho}.
	\label{eq:rel2}
\end{equation}
Inserting relations \eqref{eq:rel1} and \eqref{eq:rel2} in \eqref{eq:EL.s} 
equation \eqref{eq:Acoustic.rho} can be derived from the Euler-Lagrange equations. 

\paragraph{The system of nonlinear Maxwell equations in absence of charges.}

The nonlinear Maxwell equations in the absence of charges are of the form
\begin{align}
	\frac{\partial \B}{\partial t} + \nabla \times\E  &=0, \label{eq:Maxwell.B}\\
	\frac{\partial \D}{\partial t} - \nabla \times\H &=0, \label{eq:Maxwell.D}
\end{align}
and satisfy the involutions $ \nabla\cdot \B = 0  $ and $ \nabla \cdot\D =0 $. In this case the state vector is $ \q = (\B, \D) $, while $ \p=(\H,\E) $ denotes the corresponding dual variables.  For a generic total energy density $ \en=\en(\B,\D) $, the system \eqref{eq:Maxwell.B}--\eqref{eq:Maxwell.D} satisfies an additional conservation law given by
\begin{equation}
	\frac{\pd \en}{\pd t} + \nabla \cdot \left( \E \times \H \right) = 0.
\end{equation}
This can be easily proven by multiplying the system \eqref{eq:Maxwell.B}--\eqref{eq:Maxwell.D} with the dual variables $ \p $ and summing up all equations.
As for the classical Maxwell equations, the system can be derived from a variational principle introducing a vector potential $ \A $. The Lagrangian $ \Lambda(\A, \pd_t \A,  \nabla\A ) $ can then be reparametrized in the following variables 
\begin{equation}
	\E = -\partial_t \A, \qquad \B =  \nabla \times \A. \label{eqn.dynamic} 
\end{equation}
i.e.
\begin{equation}
	\mathcal{L}(\E,\B) =  \Lambda(\A, \pd_t \A, \nabla \A ).
\end{equation}
Applying the time derivative to the first equation in \eqref{eqn.dynamic} and then using the second relation directly leads to the first Maxwell equation, as a mere consequence of the definitions given in \eqref{eqn.dynamic}. The second one can be obtained from the Euler-Lagrange equation. Indeed, considering the relations between the derivatives of the parametrizations $ \mathcal{L} $ and $ \Lambda $, i.e.
\begin{alignat}{2}
	\frac{\pd \Lambda}{\pd \partial_t A_{k}} & =-\frac{\pd \mathcal{L}}{\pd E_k},
	\qquad
	\frac{\pd\Lambda}{\pd \partial_j A_{k}}  & = \epsilon_{jki} \frac{\pd \mathcal{L}}{\pd B_i} 
\end{alignat}
and assuming $\pd \Lambda/\pd A_k = 0$, the variation with respect to $\A$ gives the following Euler-Lagrange equations
\begin{equation}
	\frac{\pd \mathcal{L}_{E_k}}{\pd t}  + \epsilon_{kij} \pd_i \mathcal{L}_{B_j} = 0.
	\label{eq:EL.E} 
\end{equation}
Then introducing 
\begin{equation}
	\D := \mathcal{L}_{E}, \qquad
	\en(\D,\B) := \E \mathcal{L}_{\E} - \mathcal{L},
\end{equation}
as well as the following relations between the old and new parametrization
\begin{equation}
	\H = \en_{\B} = -\mathcal{L}_{\B}, 
	\qquad
	\E= \en_{\D}
\end{equation}
the second nonlinear Maxwell equation \eqref{eq:Maxwell.D} is derived from \eqref{eq:EL.E}.

\paragraph{The nonlinear Maxwell-GLM system.}
In addition, we also consider a nonlinear extension of the Maxwell-GLM system of Munz \textit{et al.} \cite{MunzCleaning} generalized to an arbitrary Lagrangian and introduced in \cite{dumbser2024variational}. The underling system, which consists of the nonlinear Maxwell equations without charges and augmented by two acoustic subsystems, has the form 
\begin{align}
	\frac{\pd \B}{\pd t} + \nabla\times\E + \nabla \xi &= 0,\label{eq:SHTC.B}\\
	\frac{\pd \phi}{\pd t} + \nabla \cdot\H &= 0, \label{eq:SHTC.phi}\\
	\frac{\pd \D}{\pd t} -\nabla\times\H + \nabla \eta  &= 0,\label{eq:SHTC.D}\\
	\frac{\pd \psi}{\pd t} + \nabla\cdot\E &= 0. \label{eq:SHTC.psi}
\end{align}
In \eqref{eq:SHTC.B}-\eqref{eq:SHTC.psi} we have introduced the state vector $\q = \left(  \B, \phi,\D, \psi \right)$, as well as the associated vector of thermodynamic dual variables $\p = \partial \mathcal{E} / \partial \q= \left(\H, \xi,\E,\eta\right)$. 
The system \eqref{eq:SHTC.B}-\eqref{eq:SHTC.psi} admits an additional conservation law for the total energy density $ \en $ of generic form.  
Indeed, multiplying \eqref{eq:SHTC.B}-\eqref{eq:SHTC.psi} with the dual variable $ \p $ and summing up 
yields the following total energy conservation law
\begin{equation}
	\frac{\partial \en}{\partial t} + \nabla \cdot \left( \E \times \H + (\eta \E + \xi \H) \right) = 0\,. 
	\label{eqn.energy} 
\end{equation}
As shown recently in \cite{dumbser2024variational}, the system of equations \eqref{eq:SHTC.B}--\eqref{eq:SHTC.psi} can be derived from a variational principle by combining the procedures briefly introduced and applied to acoustic and Maxwell systems. In particular, the state variables must be expressed in term of a vector potential $ \A $ and a scalar potential $ Z $ as
\begin{alignat}{2}
	&\E = -\pd_t \A , 						&& \varphi =-\pd_t Z\, ,\label{eqn.definitions1}\\
	&B_k = \epsilon_{kij} \pd_i A_j + \pd_k Z, \qquad 	&& \psi = \pd_k A_k\, ,\label{eqn.definitions2}
\end{alignat}
and the Lagrangian $ \Lambda(\A,Z, \pd_t \A, \pd_t Z, \nabla \A, \nabla Z) $ parameterized into $ \mathcal{L}(\E,\B,\varphi,\psi) $, defining the following relations
\begin{alignat}{2}
	\frac{\pd \Lambda}{\pd_t A_{k}} & =-\frac{\pd \mathcal{L}}{\pd E_k},
	\qquad
	&\frac{\pd\Lambda}{\pd_t Z} 		 = -\frac{\pd \mathcal{L}}{\pd \varphi},\\
	\frac{\pd\Lambda}{\pd_j A_{k}}  & = \epsilon_{jki} \frac{\pd \mathcal{L}}{\pd B_i} + \delta_{kj} \frac{\pd \mathcal{L}}{\pd \psi},
	\qquad
	&\frac{\pd\Lambda}{\pd_k Z} 		 = \frac{\pd \mathcal{L}}{\pd B_k}.
\end{alignat}
Equations \eqref{eq:SHTC.B} and \eqref{eq:SHTC.psi} are obtained as direct consequences of the definitions \eqref{eqn.definitions1} and \eqref{eqn.definitions2} by differentiating the third and fourth relations with respect to time. On the other hand, equations \eqref{eq:SHTC.D} and \eqref{eq:SHTC.phi} are Euler-Lagrange equations. 
Indeed, assuming $\pd \Lambda/\pd A_k = 0$ and $\pd \Lambda/ \pd Z = 0$, the
variation with respect to $A_k$ and $Z$ gives the Euler-Lagrange equations 
\begin{equation}
	\frac{\pd \mathcal{L}_{E_k}}{\pd t}  + \epsilon_{kij} \pd_i \mathcal{L}_{B_j} - \pd_k \mathcal{L}_\psi = 0, \qquad 
	\frac{\pd \mathcal{L}_\varphi}{\pd t}  -\pd_k \mathcal{L}_{B_k} = 0.
	\label{eq:EL.E.phi} 
\end{equation}
After introducing 
\begin{equation}
	\D := \mathcal{L}_{\E}, \qquad \phi := \mathcal{L}_\varphi,
	\qquad
	\en(\D,\B,\phi,\psi) := \E \mathcal{L}_{\E} + \varphi \mathcal{L}_\varphi - \mathcal{L},
\end{equation}
with the following relations 
\begin{equation}
	\H = \en_{\B} = -\mathcal{L}_{\B}, 
	\qquad 
	\eta = \en_{\psi} = -\mathcal{L}_{\psi},
	\qquad
	\E = \en_{\D},
	\qquad
	\xi = \en_{\phi},
\end{equation}
equations \eqref{eq:EL.E.phi} can be rewritten to derive \eqref{eq:SHTC.D} and \eqref{eq:SHTC.phi}.

\vspace{1cm}

Moreover, after a change of variables, all systems considered above can be immediately symmetrized in the dual variables and the generating potential, rewriting them in the following form
\begin{equation} 
	L_{\p\p} \frac{\partial \p}{\partial t} + H_k \frac{\partial \p}{\partial x_k} = 0.
	\label{eq:sysP}
\end{equation} 
The $ H_k $ matrices are constant in all cases considered above and read 
\begin{equation}\label{eqn.H.1}
	H_1 = \left(  \begin{matrix}
		0 & 0 & 0 & 1  \\ 
		0 & 0 & 0 & 0  \\
		0 & 0 & 0 & 0  \\
		1 & 0 & 0 & 0 \\ 
	\end{matrix}
	\right),   	 
\qquad
	H_2 = \left( \begin{array}{cccccccc} 
		0 & 0 &  0 & 0  \\ 
		0 & 0 & 0 & 1\\
		0 & 0 & 0 & 0 \\
		0 & 1 & 0 & 0 \\ 
	\end{array} \right),   	 
\qquad
	H_3 = \left( \begin{array}{cccccccc} 
		0 & 0 & 0 & 0  \\ 
		0 & 0 & 0 & 0  \\
		0 & 0 & 0 & 1  \\
		0 & 0 & 1 & 0  \\ 
	\end{array} \right),
\end{equation}
when the nonlinear acoustic system \eqref{eq:Acoustic.v}--\eqref{eq:Acoustic.rho} is considered. On the other hand for the nonlinear Maxwell equations \eqref{eq:Maxwell.B}--\eqref{eq:Maxwell.D} they are given by 
\begin{equation}\label{eqn.H.1}
	H_1 = \left(  \begin{matrix}
		0 & 0 & 0 & 0 & 0 & 0 \\ 
		0 & 0 & 0 & 0 & 0 & -1  \\
		0 & 0 & 0 & 0 & 1 &  0 \\
		0 & 0 & 0 & 0 & 0 & 0  \\ 
		0 & 0 & 1 & 0 & 0 & 0 \\
		0 & -1 & 0 & 0 & 0 & 0 \\
	\end{matrix}
	\right),   	 
\quad
	H_2 = \left( \begin{array}{cccccccc} 
		0 & 0 & 0 & 0 & 0 & 1 \\ 
		0 & 0 & 0 & 0 & 0 & 0  \\
		0 & 0 & 0 & -1 & 0 &  0 \\
		0 & 0 & -1 & 0 & 0 & 0  \\ 
		0 & 0 & 0 & 0 & 0 & 0 \\
		1 & 0 & 0 & 0 & 0 & 0 \\
	\end{array} \right),   	 
\quad
	H_3 = \left( \begin{array}{cccccccc} 
	0 & 0 & 0 & 0 & -1 & 0 \\ 
	0 & 0 & 0 & 1 & 0 & 0  \\
	0 & 0 & 0 & 0 & 0 &  0 \\
	0 & 1 & 0 & 0 & 0 & 0  \\ 
	-1 & 0 & 0 & 0 & 0 & 0 \\
	0 & 0 & 0 & 0 & 0 & 0 \\
	\end{array} \right). 
\end{equation}
Last, for the Maxwell-GLM system \eqref{eq:SHTC.B}--\eqref{eq:SHTC.psi}, the matrices $ H_k $ are as follows
\begin{equation}\label{eqn.H.1}
	H_1 = \left(  \begin{matrix}
		0 & 0 & 0 & 1 & 0 & 0 & 0 &  0 \\ 
		0 & 0 & 0 & 0 & 0 & 0 & -1 & 0 \\
		0 & 0 & 0 & 0 & 0 &  1 & 0 & 0 \\
		1 & 0 & 0 & 0 & 0 & 0 & 0 & 0 \\ 
		0 & 0 & 0 & 0 & 0 & 0 & 0 & 1 \\ 
		0 & 0 & 1 & 0 & 0 & 0 &  0 & 0 \\
		0 & -1 & 0 & 0 & 0 & 0 & 0 & 0 \\
		0 & 0 & 0 & 0 &  1 & 0 & 0 & 0 \\ 
	\end{matrix}
	\right),   	 
	\qquad
	H_2 = \left( \begin{array}{cccccccc} 
		0 & 0 &  0 & 0 & 0 & 0 & 1 & 0 \\ 
		0 & 0 & 0 & 1 & 0 & 0 & 0 &  0 \\
		0 & 0 & 0 & 0 & -1 & 0 & 0 & 0 \\
		0 & 1 & 0 & 0 & 0 &  0 & 0 & 0 \\ 
		0 & 0 & -1 & 0 & 0 & 0 & 0 & 0 \\ 
		0 & 0 & 0 &  0 & 0 & 0 & 0 & 1 \\
		1 & 0 & 0 & 0 &  0 & 0 & 0 & 0 \\
		0 & 0 & 0 & 0 & 0 & 1 & 0 & 0 \\ 
	\end{array} \right),   	 \nonumber
\end{equation}
\begin{equation}
	H_3 = \left( \begin{array}{cccccccc} 
		0 & 0 & 0 & 0 & 0 & -1 & 0 & 0 \\ 
		0 & 0 & 0 & 0 & 1 & 0 & 0 & 0 \\
		0 & 0 & 0 & 1 & 0 & 0 & 0 & 0 \\
		0 & 0 & 1 & 0 & 0 & 0 &  0 & 0 \\ 
		0 & 1 & 0 & 0 & 0 &  0 & 0 & 0 \\ 
		-1 & 0 & 0 & 0 & 0 & 0 & 0 & 0 \\
		0 & 0 & 0 &  0 & 0 & 0 & 0 & 1 \\
		0 & 0 & 0 & 0 & 0 & 0 & 1 & 0 \\ 
	\end{array} \right). 
\end{equation}
Therefore, the matrices $H_k$ are obviously symmetric in all considered cases and if the matrix $L_{\p\p}$ is symmetric and strictly positive definite, which is an immediate consequence of the assumed strict convexity of the total energy potential, we can then conclude that system \eqref{eq:sysP} is symmetric hyperbolic in the sense of Friedrichs \cite{FriedrichsSymm}.

\section{Thermodynamically compatible numerical discretizations}
\label{sec.scheme}

\subsection{Semi-discrete energy-conserving finite volume scheme on collocated meshes}
Based on the ideas outlined in \cite{HTCGPR,HTCMHD,HTCAbgrall,HTCTwoFluid,dumbser2024variational} in the following we show how to achieve the thermodynamic compatibility property \eqref{eqn.flux.comp} exactly also at the discrete level. 

We first introduce a computational domain  $\Omega \subset \mathbb{R}^d$ in $d$ space dimensions which is paved by a collocated Cartesian mesh. The common edge / face of two neighboring elements $\Omega^{\ell}$ and $\Omega^{\err}$ is denoted by $\partial \Omega^{\ell \err} = \Omega^{\ell} \cap \Omega^{\err}$ and $\mathcal{N}_\ell$ is the set of neighbors of the element $\Omega^{\ell}$. 
A semi-discrete finite volume scheme for \eqref{eqn.pde} reads
\begin{equation}
	\frac{\partial \q^{\ell}}{\partial t}  = -  \sum_{\Omega^\err \in \mathcal{N}_{\ell}}
	\frac{\left|\partial\Omega^{\ell\err}\right|}{\left|\Omega^{\ell}\right|}  \f^{\ell\err}. 
	\label{eqn.flux2dfv}
\end{equation}
To guarantee the thermodynamic compatibility between \eqref{eqn.pde} and \eqref{eq:conE} also at the discrete level, the numerical flux in the normal direction, $  \f^{\ell\err} $, must fulfill an analogous discrete compatibility condition given by
\begin{equation}
	\p^{\ell} \cdot \left( \f^{\ell\err} - \f_k^{\ell} n_k^{\ell\err} \right) +
	\p^{\err} \cdot \left( \f_k^{\err} n_k^{\ell\err} - \f^{\ell\err} \right) =  \left(F_k^{\err}  - F_k^{\ell} \right) n_k^{\ell\err}  ,\label{eqn.compatibility_2dFV}
\end{equation} 
with $\mathbf{n}^{\ell \err} = \left\{ n_k^{\ell \err} \right\}$ the unit normal vector pointing from element $\Omega^{\ell}$ to its neighbor $\Omega^{\err}$. 
Therefore, the following identity is obviously true: $\mathbf{n}^{\err\ell} = -\mathbf{n}^{\ell \err}$.
A numerical scheme of type \eqref{eqn.flux2dfv} which satisfies the property \eqref{eqn.compatibility_2dFV} preserves the total energy, see \cite{dumbser2024variational}.
Throughout this paper the scheme \eqref{eqn.flux2dfv} is integrated in time via a high order Runge-Kutta method of suitable order.

\paragraph{Compatible Abgrall-type scheme} 

As in \cite{Abgrall2018,HTCAbgrall,HTCTwoFluid,dumbser2024variational}, we employ the thermodynamically compatible Abgrall flux, which reads  
\begin{equation}
	\f^{\ell\err} = 
	\halb \left( \mathbf{f}_k^{\ell} + \mathbf{f}_k^{\err} \right) n_k^{\ell \err}  
	- \alpha^{\ell \err}  	 \left( \p^{\err} - \p^{\ell} \right).	
	\label{eqn.abgrallflux} 
\end{equation}
The scalar correction factor $\alpha^{\ell \err}$ is obtained directly by imposing the discrete compatibility condition \eqref{eqn.compatibility_2dFV} on the numerical flux and therefore is given by 
\begin{equation}
	\alpha^{\ell \err} =  \frac{ \left( F^{\err}_k - F^{\ell}_k \right) n_k^{\ell \err} + 
		\halb \left(  \p^{\err} + \p^{\ell}  \right) \cdot   \left( \mathbf{f}_k^{\ell} - \mathbf{f}_k^{\err} \right) n_k^{\ell \err}  }{ \left(  \p^{\err} - \p^{\ell}  \right)^2 }.  
	\label{eqn.alpha2d}  		
\end{equation}	
It then follows trivially from the construction that the Abgrall flux \eqref{eqn.abgrallflux} with the correction factor \eqref{eqn.alpha2d} satisfies the discrete compatibility condition \eqref{eqn.compatibility_2dFV}.  

\subsection{Fully-discrete structure-preserving staggered semi-implicit scheme}
In this section we consider a fully-discrete, energy-conserving semi-implicit scheme, which employs a staggered mesh with mimetic discrete differential operators that preserve the essential vector calculus identities exactly at the discrete level. As in \cite{dumbser2024variational}, the computational domain $ \Omega $ is covered by a Cartesian mesh composed of primal control volumes $ \Omega_c $ with cell centers $ \x_c $ and vertices $ \x_p $. On the other hand, a dual mesh can be simultaneously built whose control volumes are denoted by $ \Omega_p $ with cell centers $ \x_p $ and vertices $ \x_c $.
As reported in \cite{dumbser2024variational}, two mimetic nabla operators are defined
\begin{equation}
	\nabla_c^p \circ = \frac{1}{|\Omega_c|} \sum_{p \in \Omega_c} l_{pc} \mathbf{n}_{pc}  \circ,
	\label{eqn.nablacp}  
\end{equation} 	
and its dual
\begin{equation}
	\nabla_p^c \circ = \frac{1}{|\Omega_p|} \sum_{c \in \Omega_p} l_{pc} \mathbf{n}_{cp} \circ,
	\label{eqn.nablapc}  
\end{equation} 	
with $\mathbf{n}_{cp} = -\mathbf{n}_{pc}$. On a Cartesian grid in two space dimensions the corner normals simply read $l_{pc} \mathbf{n}_{pc} = \halb \left( \pm \Delta y, \pm \Delta x \right)^T$. We note that, due to the Gauss theorem, the following identities are obviously true 
\begin{equation}
	\sum_{p \in \Omega_c} l_{pc} \mathbf{n}_{pc} = 0, \qquad 
	\sum_{c \in \Omega_p} l_{pc} \mathbf{n}_{cp} = 0. 
	\label{eqn.gauss.disc} 
\end{equation} 
These discrete nabla operators, in the forms given by equations  \eqref{eqn.nablacp} and \eqref{eqn.nablapc} have not yet been applied to any fields. Let $\phi_p=\phi(\x_p)$ and $\A_p=\A(\x_p)$ be a scalar and vector fields, respectively, defined on the vertices of the primal grid and let $\phi_c=\phi(\x_c)$ and $\A_c=\A(\x_c)$ be their counterparts defined on the centers of the primal grid (the vertices of the dual grid), then the following discrete gradient, divergence and curl operators and their duals are naturally defined via \eqref{eqn.nablacp} and \eqref{eqn.nablapc}: 
\begin{alignat}{2}
	&\nabla_c^p \phi_p   = \frac{1}{|\Omega_c|} \sum_{p \in \Omega_c} l_{pc} \mathbf{n}_{pc} \,     \phi_p,  \qquad 
	&&\nabla_p^c \phi_c      = \frac{1}{|\Omega_p|} \sum_{c \in \Omega_p} l_{pc} \mathbf{n}_{cp} \,     \phi_c, 
	\\  	
	&\nabla_c^p \cdot  \A_p = \frac{1}{|\Omega_c|} \sum_{p \in \Omega_c} l_{pc} \mathbf{n}_{pc} \cdot  \A_p,  \qquad
	&&\nabla_p^c \cdot  \A_c = \frac{1}{|\Omega_p|} \sum_{c \in \Omega_p} l_{pc} \mathbf{n}_{cp} \cdot  \A_c, 
	\\ 
	&\nabla_c^p \times \A_p = \frac{1}{|\Omega_c|} \sum_{p \in \Omega_c} l_{pc} \mathbf{n}_{pc} \times \A_p,  \qquad
	&&\nabla_p^c \times \A_c = \frac{1}{|\Omega_p|} \sum_{c \in \Omega_p} l_{pc} \mathbf{n}_{cp} \times \A_c. 
\end{alignat} 
In addition, the following discrete vector calculus identities are obtained from some calculations. These identities form the fundamental basis of the proposed numerical scheme:
The discrete curl applied to the discrete gradient of a discrete scalar field $ \phi  $ vanishes identically, i.e. 
\begin{equation}
	\nabla_c^p \times \nabla_p^a \, \phi_a = 0, \qquad 
	\nabla_p^c \times \nabla_c^q \, \phi_q = 0,  
	\label{eqn.rot.grad} 
\end{equation}
which is the discrete analogue of the continuous identity $\nabla \times \nabla \phi = 0$.
The discrete divergence applied to the discrete curl of a vector field $ \A $ vanishes, we have 
\begin{equation}
	\nabla_c^p \cdot \nabla_p^a \times \A_a = 0, \qquad 
	\nabla_p^c \cdot \nabla_c^q \times \A_q = 0,  
	\label{eqn.div.rot} 
\end{equation}
which reflects the continuous identity $\nabla \cdot \nabla \times \mathbf{A} = 0 $ at the discrete level.
We underline that repeated indices denote summation over cells and vertices, respectively, according to the usual Einstein summation convention. The proof can be found, for example, in \cite{SIGPR} and \cite{dumbser2024variational}. 

The quantities are arranged in different and appropriate locations on the staggered grid. Generally, the state variables whose governing equation is derived from an Euler-Lagrangian equation as well as their duals are defined at the vertices of the primal mesh, which can also be seen as the barycenters of dual cells. Conversely, the state variables whose evolution equation is obtained as direct consequence of definitions and their associated dual variables are located at the centers of the primal mesh. In particular, to describe the acoustic system \eqref{eq:Acoustic.v}--\eqref{eq:Acoustic.rho} we define $ \vv^n_c=\vv(\x_c,t^n) $, $ \rho^n_p=\rho(\x_p,t^n) $ and $ s^n_p=s(\x_p,t^n) $. $ \B^n_c=\B(\x_c,t^n) $, $ \D^n_p=\D(\x_p,t^n) $ and $ \H^n_c=\H(\x_c,t^n) $ $ \E^n_p=\E(\x_p,t^n) $ denote the discrete magnetic and electric fields together with their dual variables. Finally, the scalars $ \phi $ and $ \psi $ in \eqref{eq:SHTC.phi} and \eqref{eq:SHTC.psi} and their duals $ \xi $ and $ \eta $ are represented by $ \phi^n_c=\phi(\x_c,t^n) $, $ \xi^n_c=\xi(\x_c,t^n) $ and $ \psi^n_p=\psi(\x_p,t^n) $, $ \eta^n_p=\eta(\x_p,t^n) $, respectively.


Differently from what was proposed in \cite{dumbser2024variational}, we discretize the system written in terms of the dual variables in symmetrised form \eqref{eq:sysP}.
We therefore propose the following compatible discretization of \eqref{eq:sysP}:
\begin{equation}
	\tL_{\p\p,i} \cdot \frac{\p_i^{n+1}-\p_i^n}{\Delta t} + H_h\nabla_h \tp_i =0,
	\label{eq:sysPdiscr}
\end{equation}
where 
\begin{equation}
	\tp_i = \int_{0}^{1} \p(\pi) ds; \qquad
	\left(\tL_{\p\p,i}\right)^{-1} = \int_{0}^{1} \en_{\q\q}(\pi) ds
	\label{eq:paths}
\end{equation}
with $ \pi = \q_i^n + s(\q_i^{n+1}-\q_i^n) $ the simple straight line segment path in $ \q $ variables and $ \en_{\q\q} $ the Hessian of the energy potential. Note that the subscript $ i $ indicates a proper location on the mesh and the discrete differential operator $ \nabla_h $ represents the two mimetic nabla operators defined in \eqref{eqn.nablacp} and \eqref{eqn.nablapc}.

Due to the presence of the two unknown vectors, $ \p_i^{n+1} $ and $ \tp_i $, in \eqref{eq:sysPdiscr}, a Picard loop with an initial guess of $ \p_i^{n+1,0} = \p_i^n $ is introduced, and the trick of rewriting $ \tp_i $ as
\begin{align}
	\tp_i &= \left(\halb \left( \p_i^{n+1}-\p_i^n\right)\right)^{m+1} +  \left(\tp_i-\halb \left( \p_i^{n+1}-\p_i^n\right)\right)^{m} \nonumber\\
	&= \halb \left(\p_i^{n+1,m+1}-\p_i^{n+1,m}\right)+\tp_i^m
	\label{eq:ptilde}
\end{align}
is employed.
Then, at each Picard iteration we get the following system 
\begin{equation}
	\left(\tL_{\p\p,i}\right)^m \cdot \frac{\p_i^{n+1,m+1}-\p_i^n}{\Delta t} + H_h\nabla_h\left(\halb \left(\p_i^{n+1,m+1}-\p_i^{n+1,m}\right)+\tp_i^m\right) =0,
\end{equation}
which can be rearranged as
\begin{equation}
	\left(\tL_{\p\p,i}\right)^m\! \cdot \p_i^{n+1,m+1} +\frac{\Delta t}{2}H_h\nabla_h\p_i^{n+1,m+1} =\left(\tL_{\p\p,i}\right)^m\! \cdot \p_i^{n} +\frac{\Delta t}{2} H_h\nabla_h\p_i^{n+1,m} + \Delta tH_h \nabla_h\tp_i^m.
	\label{eq:sysPicard}
\end{equation}
The resulting linear algebraic systems arising from \eqref{eq:sysPicard} are solved iteratively via a suitable Krylov subspace method, together with the Gauss quadrature rule to calculate the path integrals in time arising in \ref{eq:paths}. Once $ \p_i^{n+1,m+1} $ is computed, $ \q_i^{n+1,m+1} $ can be derived and $ \en^{n+1,m+1} $ calculated. The iterations stop when $ |\en^{n+1,m+1} - \en^{n+1,m}| < tol $, with $ tol = 10^{-15} $.

\begin{property}
	\label{prop.energy} 	
	In the case of periodic boundaries the scheme \eqref{eq:sysPdiscr} conserves the global discrete total energy exactly, i.e. 
	$$\mathcal{E}^{n+1} = \mathcal{E}^n.$$  
\end{property} 
\begin{proof}
	We compute the dot product of $ \tp_i $ and the scheme \eqref{eq:sysPdiscr} leading to 
	\begin{equation}
		\tp_i \cdot \left(\tL_{\p\p,i} \cdot \frac{\p_i^{n+1}-\p_i^n}{\Delta t} + H_h\nabla_h \tp_i \right)=0.
		\label{eq:sysPdiscr2}
	\end{equation}
    We then observe that in order to guarantee a thermodynamically compatible scheme, the matrix $ \tL_{\p\p,i} $ must satisfy the following Roe property
    \begin{equation}
    	\tL_{\p\p,i} \cdot \frac{\p_i^{n+1}-\p_i^n}{\Delta t} = \frac{\q_i^{n+1}-\q_i^n}{\Delta t}
    \end{equation}
    and the following discrete chain rule must be fulfilled
    \begin{equation}
    	\tp_i \cdot \frac{\q_i^{n+1}-\q_i^n}{\Delta t} = \frac{\en_i^{n+1}-\en_i^n}{\Delta t}.
    \end{equation}
    The first requirement can be easily proven by considering the relations
    \begin{equation}
    	\p_i^{n+1}-\p_i^n = \int_{\p_i^n}^{\p_i^{n+1}} \, d\p = \int_{\q_i^n}^{\q_i^{n+1}} \frac{\partial \p}{\partial \q} \, d\q = \int_{\q_i^n}^{\q_i^{n+1}} \frac{\partial \en_\q}{\partial \q} \, d\q = \int_{\q_i^n}^{\q_i^{n+1}} \en_{\q\q} \, d\q.
    	\label{eq:enrel2}
    \end{equation}
    Then, we introduce the straight line segment path $ \pi = \q_i^n + s(\q_i^{n+1}-\q_i^n) $ and we recall the definition of $ \tL_{\p\p,i} $ in \eqref{eq:paths}. In this way, the last integral in \eqref{eq:enrel2} can be described as
    \begin{equation}
    	\int_{\q_i^n}^{\q_i^{n+1}} \en_{\q\q} \, d\q = \int_{0}^{1} \en_{\q\q}(\pi) \frac{\partial \pi}{\partial s} \, ds = \left(\tL_{\p\p,i} \right)^{-1}\, \left(\q_i^{n+1}-\q_i^n\right),
    	\label{eq:enrel3}
    \end{equation}
    and so from \eqref{eq:enrel2} and \eqref{eq:enrel3} we have
    \begin{equation}
    	\tL_{\p\p,i} \left(\p_i^{n+1}-\p_i^n\right) = \left(\q_i^{n+1}-\q_i^n\right).
    \end{equation}
    Similarly, the second condition can also be proven.
	We start from the relation
	\begin{equation}
		\en_i^{n+1}-\en_i^n = \int_{\en_i^n}^{\en_i^{n+1}} \, d\en = \int_{\q_i^n}^{\q_i^{n+1}} \frac{\partial \en}{\partial \q} \, d\q = \int_{\q_i^n}^{\q_i^{n+1}} \p \, d\q.
		\label{eq:enrel1}
	\end{equation}
	Then, introducing again the line segment path $ \pi = \q_i^n + s(\q_i^{n+1}-\q_i^n) $ and recalling the definition of $ \tp_i $ in \eqref{eq:paths}, leads to  
	\begin{equation}
		\int_{\q_i^n}^{\q_i^{n+1}} \p \, d\q = \int_{0}^{1} \p(\pi) \frac{\partial \pi}{\partial s} \, ds = \tp_i \, \left(\q_i^{n+1}-\q_i^n\right).
	\end{equation}
	At this point the scheme, \eqref{eq:sysPdiscr2} can be rewritten as
	\begin{equation}
		\en_i^{n+1} -\en^n_i = -\tp_i \cdot H_h\nabla_h \tp_i.
		\label{eq:diffenergy}
	\end{equation}
	We now sum up equation \eqref{eq:diffenergy} over all the elements. To prove that this sum is a null identity, we consider the systems described above individually. In particular, for the nonlinear acoustic system we get
	\begin{eqnarray} 
		 \en^{n+1} - \en^n &&= 
		- \Delta t \sum_c \tilde{\uu}_c\cdot\nabla_c^p \tilde{s}_p - \Delta t\sum_p\tilde{s}_p\nabla_p^c\cdot\tilde{\uu}_c
		\nonumber \\
		&&= - \Delta t \sum_c \tilde{\uu}_c\cdot\sum_{p \in \Omega_c} l_{pc} \mathbf{n}_{pc} \tilde{s}_p - \Delta t\sum_p\tilde{s}_p\sum_{c \in \Omega_p} l_{pc} \mathbf{n}_{cp}\cdot\tilde{\uu}_c 
		\nonumber \\
		&&= - \Delta t \sum_c \sum_{p \in \Omega_c} l_{pc} \mathbf{n}_{pc} \cdot \tilde{s}_p\tilde{\uu}_c - \Delta t\sum_p\sum_{c \in \Omega_p} l_{pc} \mathbf{n}_{cp}\cdot\tilde{s}_p\tilde{\uu}_c =0
	\end{eqnarray} 
	since $\mathbf{n}_{cp}=-\mathbf{n}_{pc}$. 
	Similarly, considering system \eqref{eq:Maxwell.B}-- \eqref{eq:Maxwell.D} yields
	\begin{eqnarray} 
		\en^{n+1} - \en^n &&= 
		- \Delta t \sum_c \tilde{\H}_c\cdot\nabla_c^p \times\tilde{\E}_p + \Delta t\sum_p\tilde{\E}_p\cdot\nabla_p^c\times\tilde{\H}_c
		\nonumber \\
		&&= - \Delta t \sum_c \tilde{\H}_c\cdot\sum_{p \in \Omega_c} l_{pc} \mathbf{n}_{pc} \times\tilde{\E}_p + \Delta t\sum_p\tilde{\E}_p\cdot\sum_{c \in \Omega_p} l_{pc} \mathbf{n}_{cp}\cdot\tilde{\H}_c 
		\nonumber \\
		&&= - \Delta t \sum_c \sum_{p \in \Omega_c} l_{pc} \mathbf{n}_{pc} \cdot \tilde{\E}_p\times\tilde{\H}_c + \Delta t\sum_p\sum_{c \in \Omega_p} l_{pc} \mathbf{n}_{cp}\cdot\tilde{\H}_p\times\tilde{\E}_c =0
	\end{eqnarray}
	since $\tilde{\B}_c \times \tilde{\E}_p = -\tilde{\E}_p \times \tilde{\B}_c$ and $\mathbf{n}_{cp}=-\mathbf{n}_{pc}$ and thus all pairwise interactions between $p$ and $c$ cancel.
	The same can be proven taking into account the Maxwell-GLM system. Indeed,
	\begin{eqnarray} 
		&& \mathcal{E}^{n+1} - \mathcal{E}^n =
		\nonumber \\  
		&=& - \Delta t \sum_c \tilde{\H}_c\cdot\nabla_c^p \times\tilde{\E}_p - \Delta t \sum_c \tilde{\H}_c\cdot\nabla_c^p\tilde{\xi}_p - \Delta t \sum_p \tilde{\xi}_p\nabla_p^c \cdot\tilde{\H}_c \nonumber\\
		&&+ \Delta t\sum_p\tilde{\E}_p\cdot\nabla_p^c\times\tilde{\H}_c - \Delta t\sum_p\tilde{\E}_p\cdot\nabla_p^c \tilde{\eta}_c - \Delta t \sum_c \tilde{\eta}_c\nabla_c^p\cdot\tilde{\E}_p
		\nonumber\\
		&=& - \Delta t \sum_c \tilde{\H}_c\cdot\sum_{p \in \Omega_c} l_{pc} \mathbf{n}_{pc}\times\tilde{\E}_p - \Delta t \sum_c \tilde{\H}_c\cdot\sum_{p \in \Omega_c} l_{pc} \mathbf{n}_{pc}\tilde{\xi}_p - \Delta t \sum_p \tilde{\xi}_p\sum_{c \in \Omega_p} l_{pc} \mathbf{n}_{cp} \cdot\tilde{\H}_c \nonumber\\
		&& + \Delta t\sum_p\tilde{\E}_p\cdot\sum_{c \in \Omega_p} l_{pc} \mathbf{n}_{cp}\times\tilde{\H}_c - \Delta t\sum_p\tilde{\E}_p\cdot\sum_{c \in \Omega_p} l_{pc} \mathbf{n}_{cp} \tilde{\eta}_c - \Delta t \sum_c \tilde{\eta}_c\sum_{p \in \Omega_c} l_{pc} \mathbf{n}_{pc}\cdot\tilde{\E}_p\nonumber\\
		&=& - \Delta t \sum_c \sum_{p \in \Omega_c} l_{pc} \mathbf{n}_{pc}\cdot\tilde{\E}_p\times\tilde{\H}_c - \Delta t \sum_c \sum_{p \in \Omega_c} l_{pc} \mathbf{n}_{pc}\cdot\tilde{\xi}_p\tilde{\H}_c - \Delta t \sum_p \sum_{c \in \Omega_p} l_{pc} \mathbf{n}_{cp} \cdot\tilde{\xi}_p\tilde{\H}_c \nonumber\\
		&& + \Delta t\sum_p\sum_{c \in \Omega_p} l_{pc} \mathbf{n}_{cp}\cdot\tilde{\H}_c\times\tilde{\E}_p- \Delta t\sum_p\sum_{c \in \Omega_p} l_{pc} \mathbf{n}_{cp}\cdot \tilde{\eta}_c\tilde{\E}_p- \Delta t \sum_c \sum_{p \in \Omega_c} l_{pc} \mathbf{n}_{pc}\cdot\tilde{\E}_p\tilde{\eta}_c \nonumber\\
		&=& 0,
	\end{eqnarray}
	since $\tilde{\B}_c \times \tilde{\E}_p = -\tilde{\E}_p \times \tilde{\B}_c$ and $\mathbf{n}_{cp}=-\mathbf{n}_{pc}$. 
	This demonstrates the total energy conservation and, consequently, the nonlinear energy stability of the proposed numerical scheme for all systems under study.
\end{proof}
 
\begin{property}
	\label{prop.maxwell} 	
	For discrete initial data that satisfy $\nabla_p^c \cdot \B_c^n = 0$, $\nabla_c^p \cdot \D_p^n = 0$ at time $t=0$, i.e. for $n=0$, we have $\nabla_p^c \cdot \B_c^{n+1} = 0$ and $\nabla_c^p \cdot \E_p^{n+1} = 0$ for all times $t>0$. 
\end{property} 
\begin{proof}
	The proposed scheme \eqref{eq:sysPdiscr2} applied to \eqref{eq:Maxwell.B}--\eqref{eq:Maxwell.D} yields
	\begin{align}
		(\tL_{\H\H})^m_c \H^{n+1,m+1}_c \! + \!\frac{\Delta t}{2}\nabla_c^p\times\E_p^{n+1,m+1} &\! =\! 
		(\tL_{\H\H})^m_c \H^{n}_c - \Delta t \nabla_c^p\times\tE_p^m \! + \frac{\Delta t}{2}\nabla_c^p\times\E_p^{n+1,m} \label{eq:Maxwell.picH}\\
		(\tL_{\E\E})^m_p \E^{n+1,m+1}_p \! - \!\frac{\Delta t}{2}\nabla_p^c\times\H_c^{n+1,m+1} & \! =\! 
		(\tL_{\E\E})^m_p \E^{n}_p + \Delta t \nabla_p^c\times\tH_c^m \!- \frac{\Delta t}{2}\nabla_p^c\times\H_c^{n+1,m}\!\!.\label{eq:Maxwell.picE}
	\end{align}
	Then, applying the discrete divergence to \eqref{eq:Maxwell.picH} and \eqref{eq:Maxwell.picE} and using the discrete identities contained in \eqref{eqn.div.rot}, one gets
	\begin{align}
		(\tL_{\H\H})^m_c \, \H^{n+1,m+1}_c  & = 
		(\tL_{\H\H})^m_c \, \H^{n}_c,\\
		(\tL_{\E\E})^m_p \, \E^{n+1,m+1}_p & =
		(\tL_{\E\E})^m_p \, \E^{n}_p .
	\end{align}
	Since $ \tL_{\p\p}  \cdot(\p^{n+1}-\p^n) = \q^{n+1}-\q^n$ and $\nabla_p^c \cdot \B_c^n = 0$, $\nabla_c^p \cdot \D_p^n = 0$ with $ n=0 $, then we obtain the desired result, i.e. $ \nabla_p^c \cdot \B_c^{n+1} = 0$, $\nabla_c^p \cdot \D_p^{n+1} = 0$. In other words, for compatible initial data, the scheme preserves the divergence-free condition of the fields exactly at the discrete level at all times. 
\end{proof}

\begin{property}
	\label{prop.acoustic} 	
	If the field $ \vv_c^n $ satisfies the constraint $ \nabla_p^c \times \vv_c^n$ at initial time, i.e. $ n=0 $, then it satisfies the differential constraint for all the times, i.e. $\nabla_p^c \times \vv_c^{n+1} = 0$. 
\end{property} 
\begin{proof}
	The proposed scheme \eqref{eq:sysPdiscr2} applied to \eqref{eq:Acoustic.v}--\eqref{eq:Acoustic.rho} yields
	\begin{align}
		(\tL_{\uu\uu})^m_c \uu^{n+1,m+1}_c + \frac{\Delta t}{2}\nabla_c^p s_p^{n+1,m+1} &= 
		(\tL_{\uu\uu})^m_c \uu^{n}_c -\Delta t \nabla_c^p\ts_p^m+ \frac{\Delta t}{2}\nabla_c^p s_p^{n+1,m},\label{eq:Acoustic.picv}\\
		(\tL_{ss})^m_ps_p^{n+1,m+1} + \frac{\Delta t}{2}\nabla_p^c\cdot\uu_c^{n+1,m+1} &=
		(\tL_{ss})^m_ps_p^{n} - \Delta t \nabla_p^c\cdot\tuu_c^m + \frac{\Delta t}{2}\nabla_p^c\cdot\uu_c^{n+1,m}.\label{eq:Acoustic.pics}
	\end{align}
	Since $ \tL_{\p\p}  \cdot(\p^{n+1}-\p^n) = \q^{n+1}-\q^n$ and after applying the discrete curl operator to \eqref{eq:Acoustic.picv}, using the discrete identities contained in \eqref{eqn.rot.grad} and given the hypothesis $ \nabla_p^c \times \vv_c^n$ for $ n=0$, leads to $ \nabla_p^c \times \vv_c^n$ for $ n>0$.
	Therefore, the scheme preserves the curl-free condition of the field $\vv$ exactly at the discrete level for all times.    
\end{proof}

\section{Numerical results}
\label{sec.results}

\subsection{Nonlinear Maxwell system: discrete total energy and divergence preservation}
\begin{figure}[!ht]
	\centering
	\includegraphics[width=0.45\textwidth]{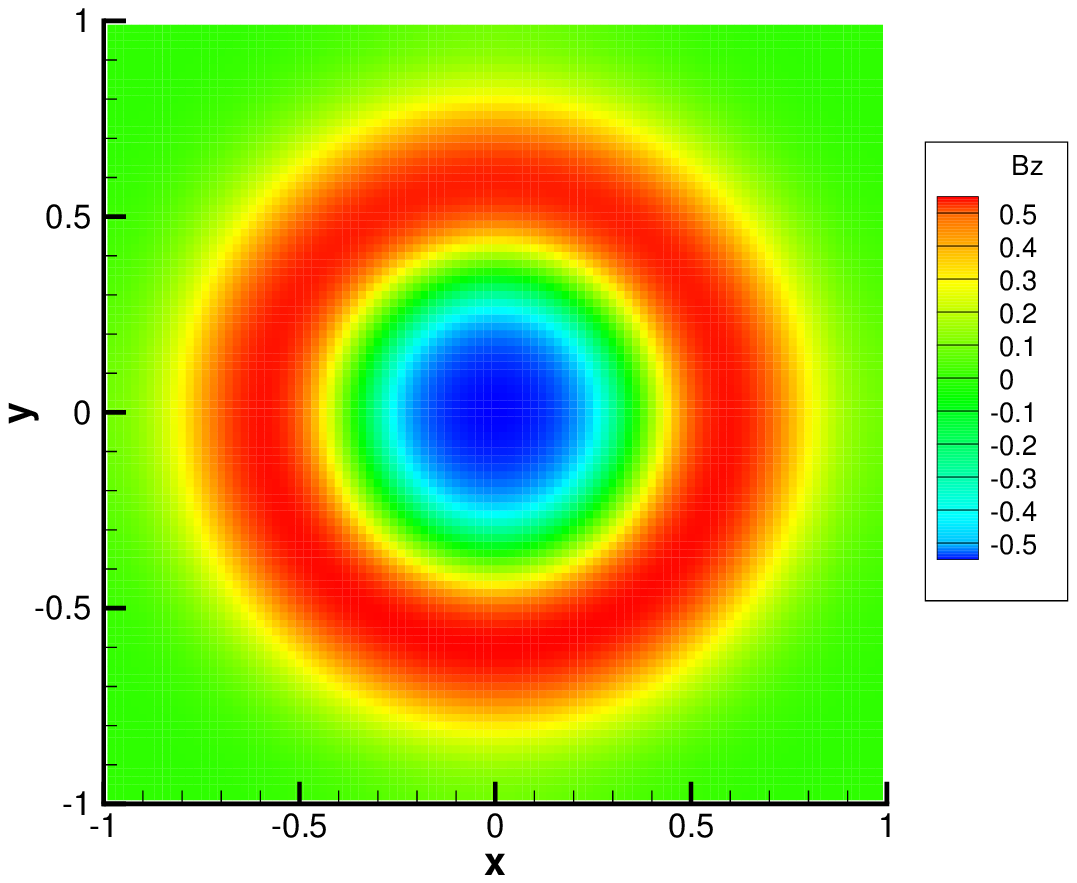}  \quad 
	\includegraphics[width=0.45\textwidth]{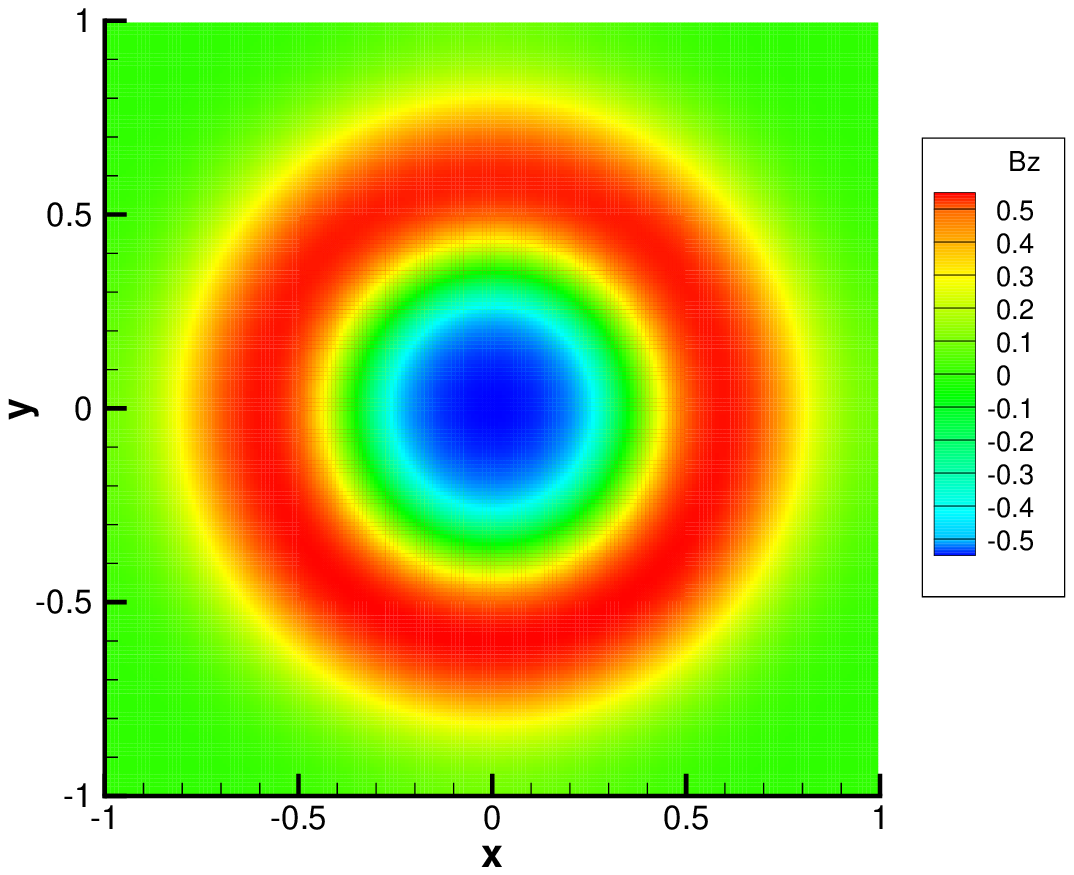}	
	\caption{The results of the nonlinear Maxwell system at $ t = 0.5 $ obtained using both the SIMM scheme (left) and HTC scheme (right), with $ 100 \times 100 $ and $ 200 \times 200 $ elements, respectively. }
	\label{fig:res.nonMaxwell}
\end{figure}

\begin{figure}[!ht]
	\centering
	\includegraphics[ width=0.45\textwidth]{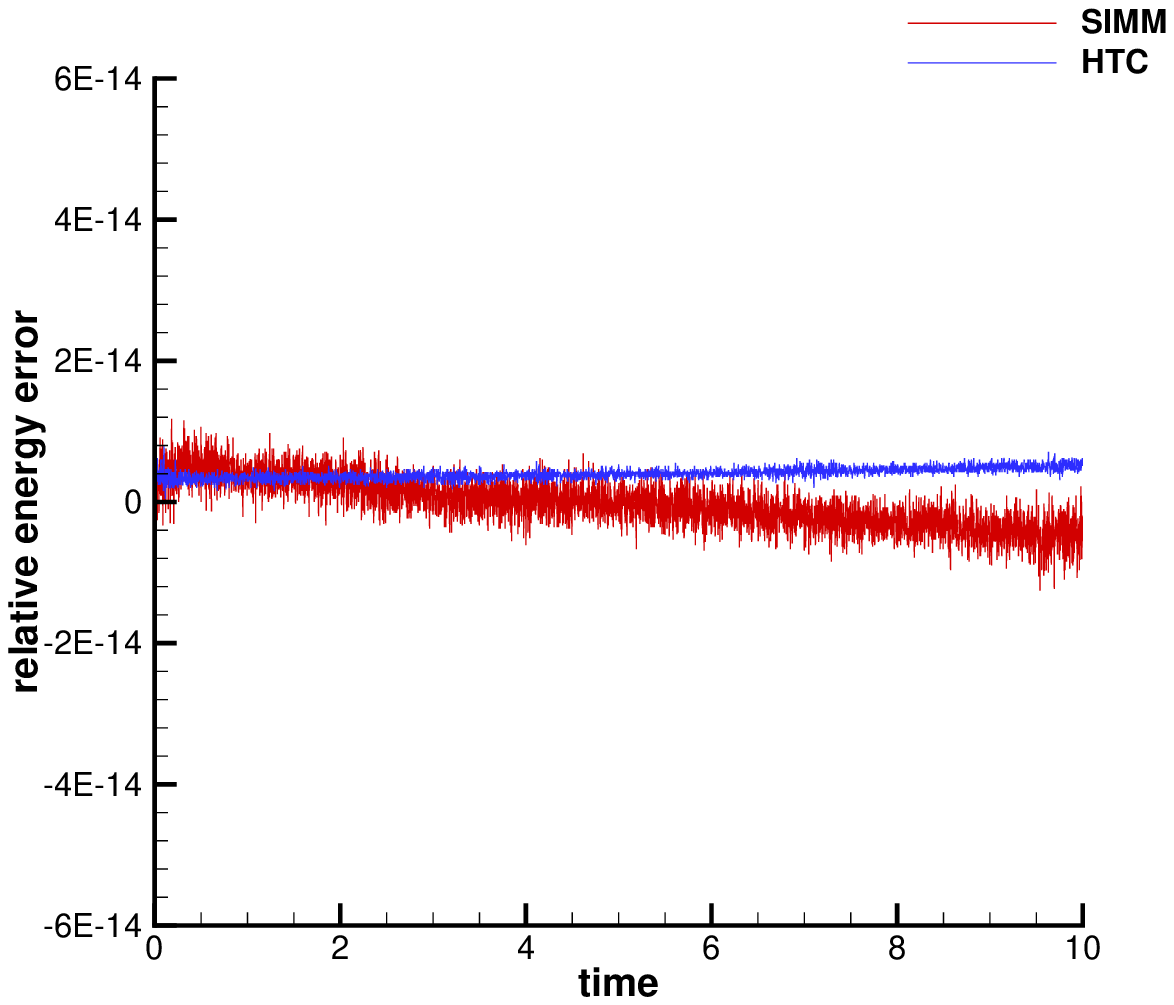}  \quad 
	\includegraphics[ width=0.45\textwidth]{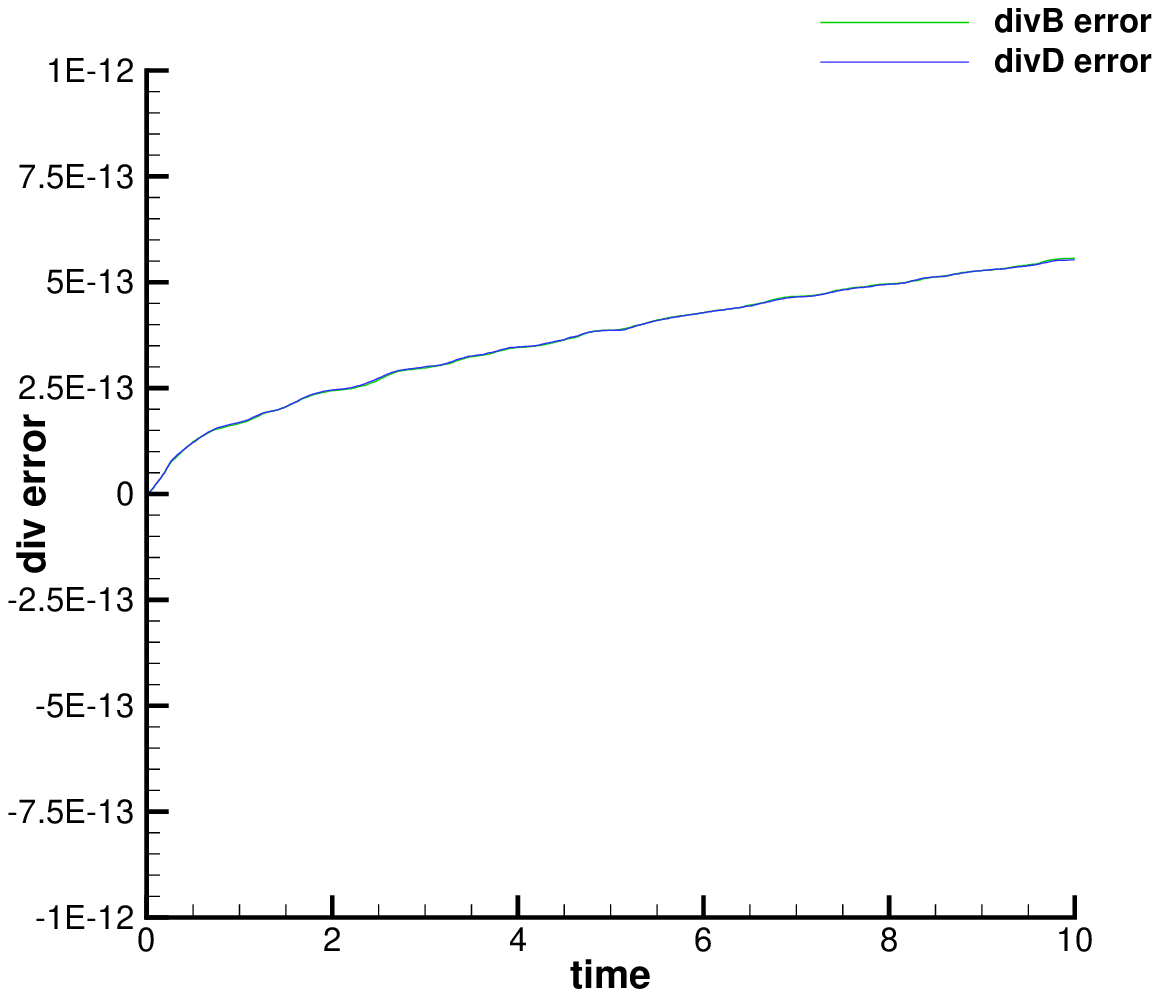}
	\caption{Left panel: Temporal evolution of the relative total energy error $\mathcal{E}^n/\mathcal{E}^0 - 1$ for initial data that is compatible with the vacuum Maxwell equations employing the semi-implicit scheme as well as the HTC scheme. Right panel: time series of the divergence errors of the magnetic and electric field for the SIMM scheme and for initial data that is compatible with the vacuum Maxwell equations. }
	\label{fig:ErrEn.nonMaxwell}
\end{figure}

As a first benchmark we consider the nonlinear Maxwell system of the form
\begin{align}
	\partial_t \B + \nabla\times\E &=0,\\
	\partial_t \D -\nabla\times\H &=0.
\end{align}
together with the following total energy density 
\begin{equation}
	\en = \B^2 + \D^2 + 0.01 \left(B_1\B^2 + D_1\D^2\right).
\end{equation}
In case of compatible initial data, i.e. $ \nabla\cdot\B =\nabla\cdot\D =0 $ at $ t=0 $, the system preserves the involutions $ \nabla\cdot\B =\nabla\cdot\D =0 $ also for all times $ t>0 $.
Both the semi-discrete HTC scheme and the fully discrete semi-implicit scheme are exactly energy
conservative. In addition, the proposed structure-preserving staggered semi-implicit scheme is also exactly compatible with the divergence free condition of the electric and the magnetic field provided the initial data are compatible. To test these properties of the scheme, we therefore solve the following test problem, given by the initial data
\begin{equation}
	\B(x, y)  = \mathbf{B_0} \exp \left( -\halb \x^2/\sigma^2 \right), \qquad  \D(x, y)  = \mathbf{D_0} \exp \left( -\halb \x^2/\sigma^2 \right).	\nonumber \\
	\label{eqn.ic.gauss}
\end{equation}
Simulations are run for all schemes in the domain $\Omega=[-1,1]^2$ until a final time of $t=10$. We choose initial data that are compatible with the original vacuum Maxwell equations. We therefore set $\mathbf{B_0} = \mathbf{D_0} = (0,0,10^{-2})$.  
The HTC scheme solves the problem on a fixed mesh composed of $ 200\times200 $ elements evolving in time the solution according to a seventh-order accurate Runge-Kutta method with a time step restricted by the CFL stability condition with CFL $ =0.5$. On the other hand, the fully-discrete semi-implicit scheme employs a $ 100\times100 $ element mesh and $ \Delta t_{\text{max}}=0.001 $.
The numerical results obtained at $ t=0.5  $ for both schemes are depicted in Figure \ref{fig:res.nonMaxwell}.
The left panel of Figure \ref{fig:ErrEn.nonMaxwell} shows the temporal evolution of the relative total energy error $\mathcal{E}^n/\mathcal{E}^0 - 1$ provided by both schemes, while in the right panel we report the time evolution of the divergence errors of the magnetic and the electric field obtained with the SIMM scheme. The method has been found to be effective in preserving total energy and involutions. 


\subsection{Nonlinear acoustic system:  discrete total energy and curl preservation}
\begin{figure}[!ht]
	\centering
	\includegraphics[width=0.45\textwidth]{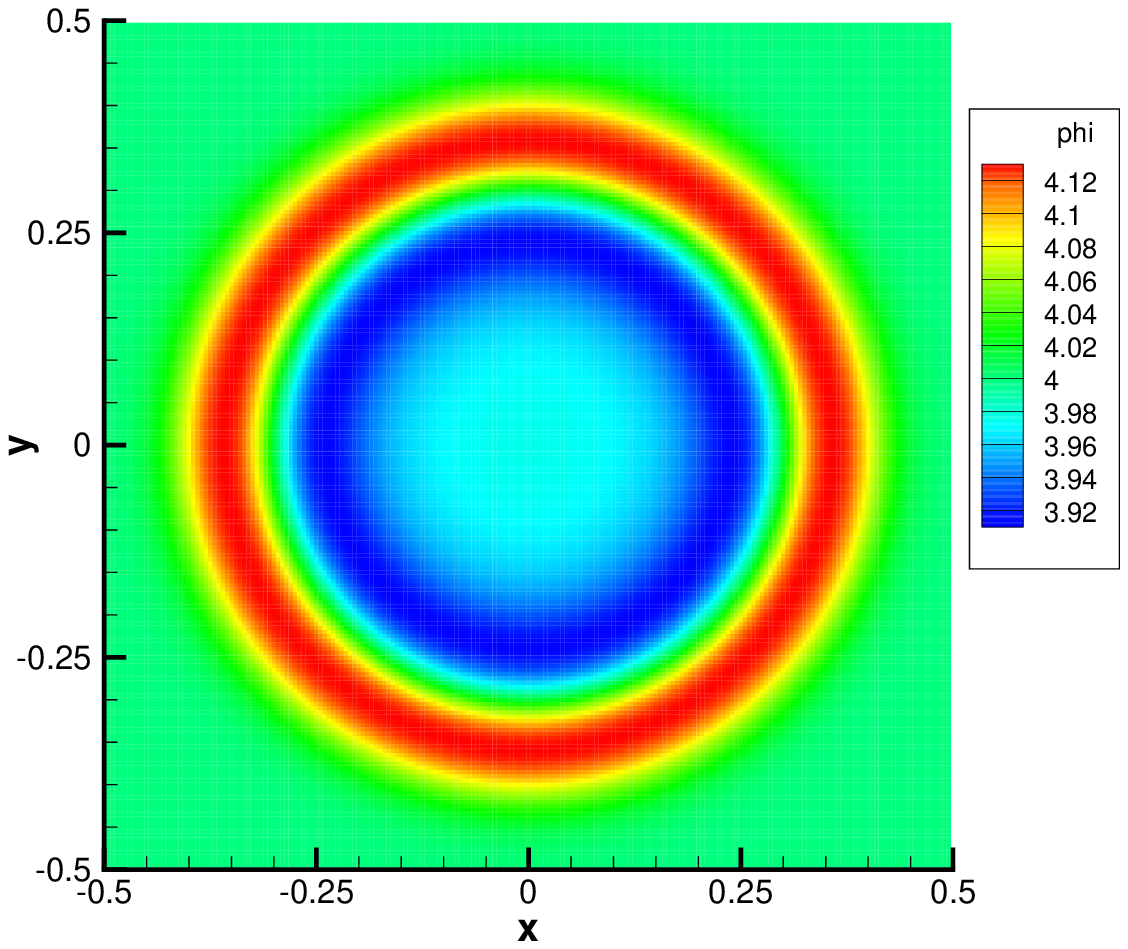}  \quad 
	\includegraphics[width=0.45\textwidth]{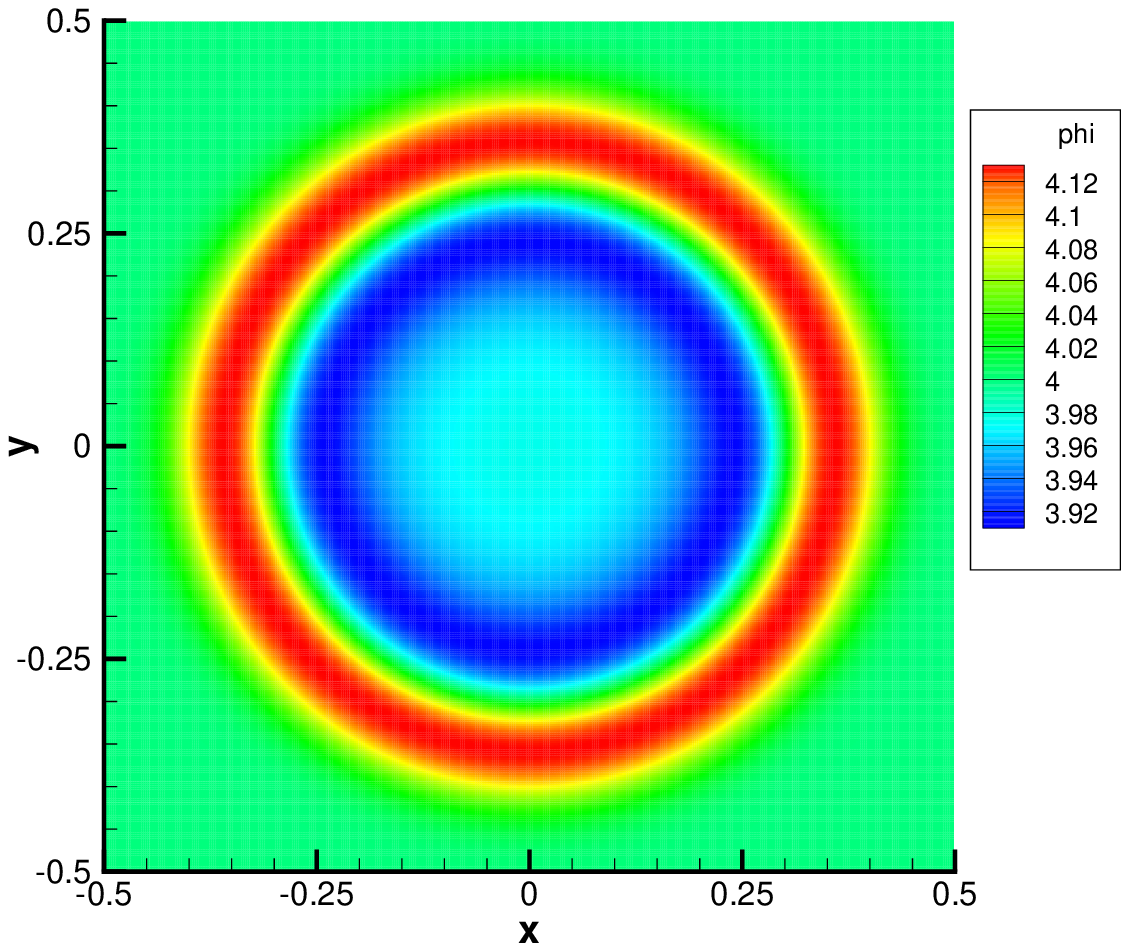}
	\caption{The results of the nonlinear acoustic subsystem at $ t = 0.25 $ obtained using both the SIMM scheme (left) and the HTC scheme (right), with $ 200 \times 200 $ and $ 400 \time 400 $ elements, respectively.}
	\label{fig:res.nonAcoustic}
\end{figure}

\begin{figure}[!ht]
	\begin{center}
	\includegraphics[width=0.45\textwidth]{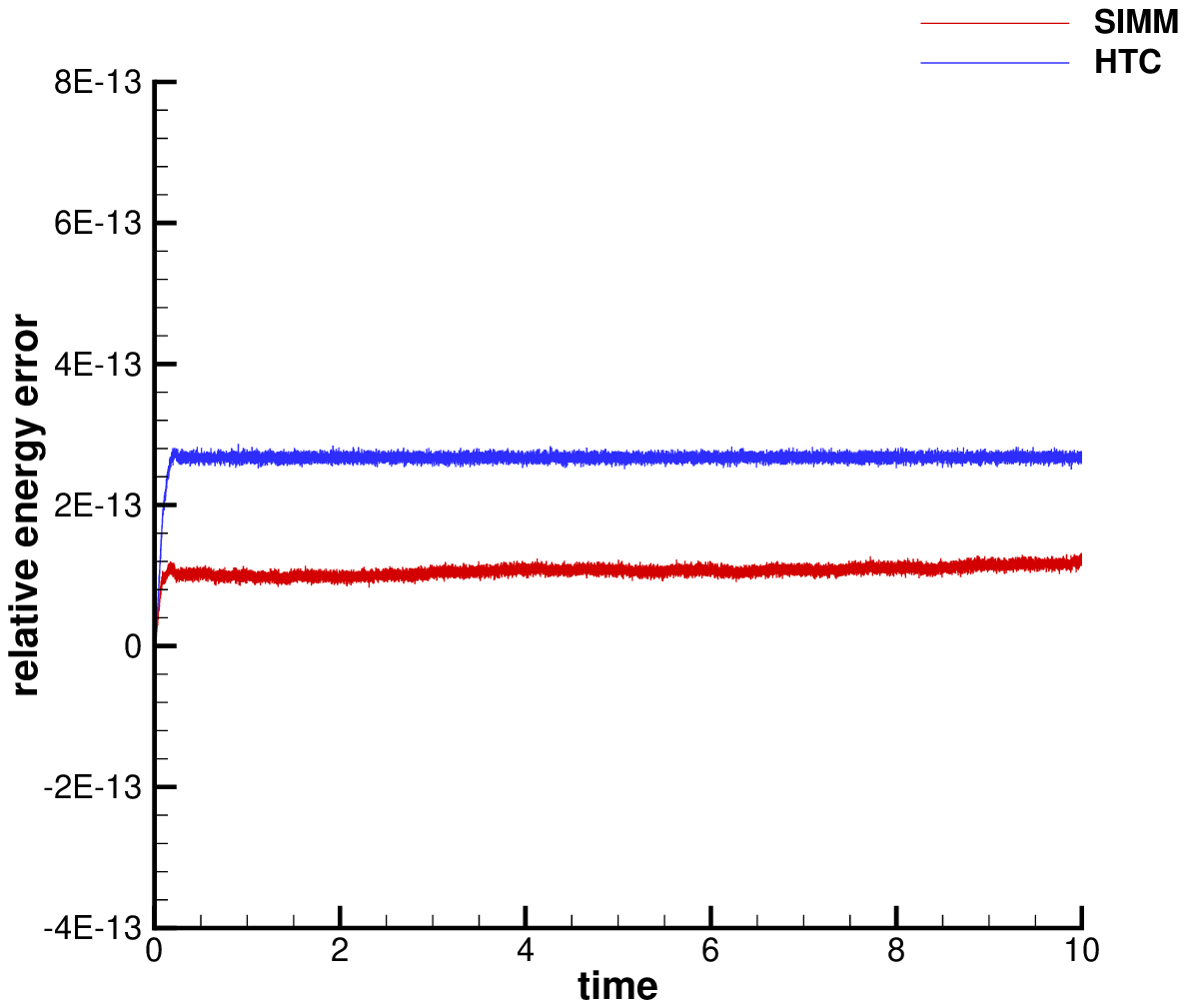}  \quad 
	\includegraphics[width=0.45\textwidth]{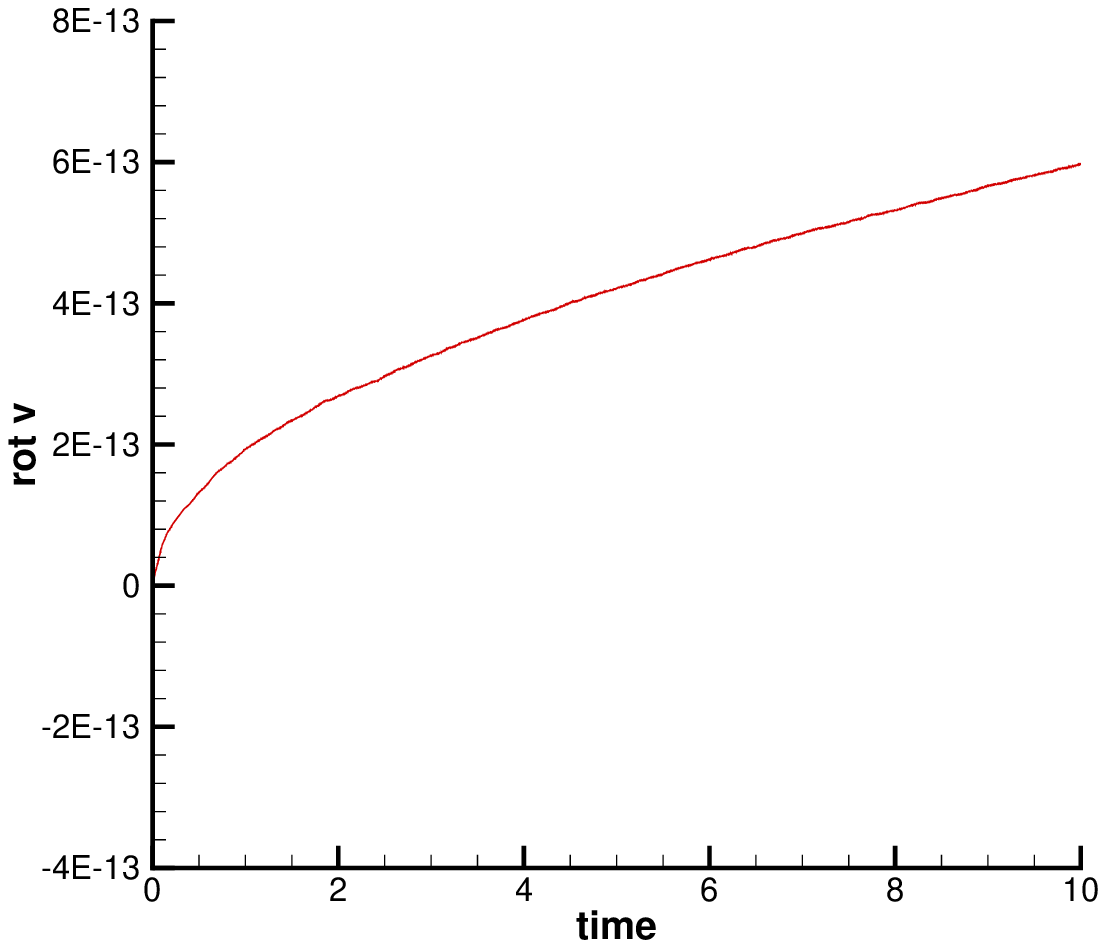}
	\end{center}
	\caption{Left panel: Temporal evolution of the relative total energy error $\mathcal{E}^n/\mathcal{E}^0 - 1$ for the nonlinear acoustic equations employing the semi-implicit scheme as well as the HTC scheme. Right panel: time series of the curl errors of the $ \vv $ field for the SIMM scheme. }
	\label{fig:ErrEn.nonAcoustic}
\end{figure}

In the second test we consider the acoustic system of the form
\begin{align}
	\partial_t \vv +\nabla s &=0,\\
	\partial_t \rho + \nabla\cdot\uu &=0.
\end{align}
The above system admits an extra conservation law for the total energy density 
\begin{equation}
	\mathcal{E} = \halb \vv^2  + \frac{ \rho^{\gamma+1}}{\gamma(\gamma+1)}
\end{equation} 
of the type
\begin{equation}
	\partial_t \en+ \nabla(\uu s) =0.
	\label{eq:en_acoustic}
\end{equation}
We recall that in this case $ s = \partial\en/\partial\rho $ and $ \uu=\partial\en/\partial\vv $.
For acoustic-compatible initial data $\nabla \times \vv=0$, the field $\vv$ preserves the stationary differential constraint  $\nabla \times \vv = 0$ for all times. This property is also satisfied by the semi-implicit scheme. 
To verify it, we consider the following initial data 
\begin{equation}
	\vv(\x, 0)  = 0, \qquad  p(\x,0) = p_0 \exp\left(-\frac{1}{2}\x^2/\sigma^2\right) + 4.
\end{equation}
The computational domain $\Omega=[-0.5,0.5]^2$ is discretized with uniform grids composed of $ 200\times 200 $ elements in case of the fully discrete semi-implicit scheme (SIMM) and of $ 400\times 400 $ elements when the semi-dicrete HTC scheme is used. Simulations are run until a final time of $t=10$. The SIMM scheme evolves the solution in time with a fixed time step of size $ \Delta t = 0.001$, while the HTC scheme employed a seventh-order Runge Kutta method where the time step is forced to obey the $\mathrm{CFL}$ stability condition with $\mathrm{CFL} =0.4 $. The parameter $ \gamma $ in \eqref{eq:en_acoustic} is chosen to be $ \gamma =1.4 $, whereas for the the initial data we set $ p_0 = 1 $ and $ \sigma = 0.05 $.
The computational results obtained for both schemes at time $ t=0.25 $ are depicted in Figure \ref{fig:res.nonAcoustic}. 
In addition, the left panel of Figure \ref{fig:ErrEn.nonAcoustic} shows the temporal evolution of the relative total energy error $\mathcal{E}^n/\mathcal{E}^0 - 1$ provided by both schemes, while the right one reports the time evolution of the curl errors of the $ \vv $ field obtained with the SIMM scheme. As expected, the method conserves total energy perfectly even over long times and the discrete velocity field remains curl-free up to machine precision, when the SIMM method is employed.

\subsection{Nonlinear Maxwell-GLM system: discrete total energy conservation}
\begin{figure}[!h]
	\centering
	\includegraphics[width=0.45\textwidth]{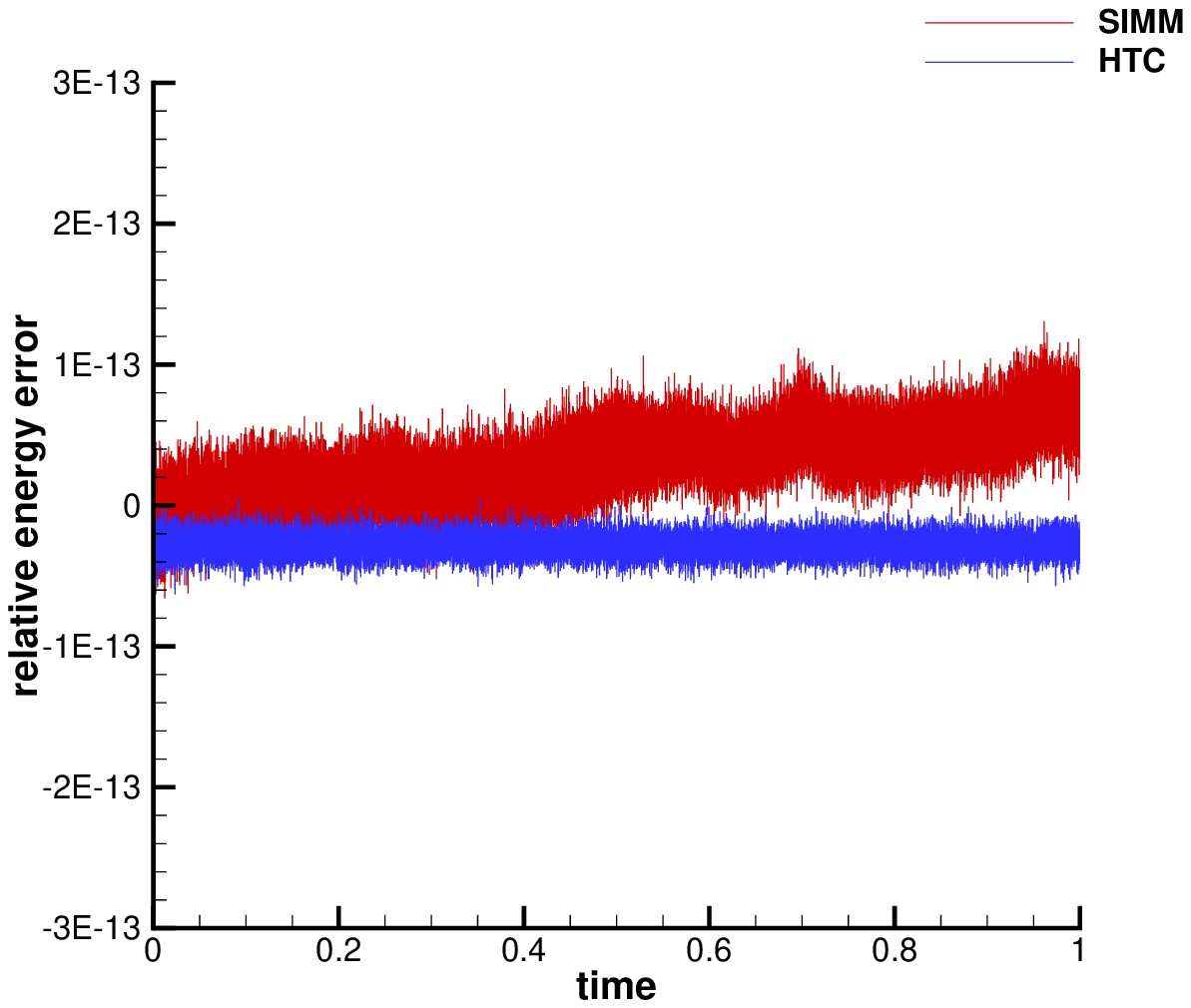}  
	\caption{Temporal evolution of the relative total energy error $\mathcal{E}^n/\mathcal{E}^0 - 1$ for the nonlinear Maxwell-GLM equations employing the semi-implicit scheme as well as the HTC scheme.  }
	\label{fig:ErrEn.MaxwellGLM}
\end{figure}

We consider the following initial data of a planar wave traveling in the direction $\normal = (1,0,0)$. The initial condition is given by 
\begin{eqnarray}
	& \B(\x, 0)  = \mathbf{B_0} \sin\left( 2\pi (x_1+x_2) \right), \qquad & \phi(\x,0) = 0.125\sin \left(2\pi (x_1+x_2)\right),	\\
	& \D(\x, 0)  = \mathbf{D_0} \sin\left( 2\pi (x_1+x_2) \right), \qquad & \psi(\x,0) = 0.25  \sin\left(2\pi (x_1+x_2)\right), 
\end{eqnarray}
with $\mathbf{B_0}=(0.125, 0, 0.5)^T $, $ \mathbf{D_0}=(0.25, 0.5, 0)^T $.
We assume the following relation between $ \B $ and its dual $ \H $
\begin{equation}
	\H = \mu_0 \B + \frac{B_1^2+B_3^2}{2}\begin{bmatrix}
		1 & 0 & 0\\
		0 & 0 & 0\\
		0 & 0 & 1
	\end{bmatrix}\B
\end{equation}
and between $ \D $ and its dual $ \E $
\begin{equation}
	\E = \mu_0 \D + \frac{D_1^2+D_3^2}{2}\begin{bmatrix}
		1 & 0 & 0\\
		0 & 0 & 0\\
		0 & 0 & 1
	\end{bmatrix}\D,
\end{equation}
where $ \mu_0 = 1 $.
In this way we obtain a total energy density given by
\begin{equation}
	\en = \halb \left( \D^2 + \B^2 + \psi^2 + \phi^2 \right) + \frac{1}{8} \left(B_1^4 + B_3^4+D_1^4 + D_3^4\right) + \frac{1}{4} \left(B_1^2 B_3^2+D_1^2 D_3^2\right).
\end{equation}
Figure \ref{fig:ErrEn.MaxwellGLM} shows the temporal evolution of the relative total energy error $\mathcal{E}^n/\mathcal{E}^0 - 1$ provided by the fully-discrete semi-implicit scheme (SIMM) on a fixed mesh composed of $2000 \times 25$ elements with $\Delta t_{\text{max}} = 0.0000125$ and by the HTC scheme on a fixed mesh composed of $ 38400 \times 50 $ elements employing a $ 7 $th-order Runge Kutta scheme. The total energy conservation is observed to be achieved at all times within the range of machine precision.

\section{Conclusions}
\label{sec.Conclusions}
In this paper, we have presented two types of structure-preserving schemes for nonlinear symmetric hyperbolic and thermodynamically compatible (SIHTC) systems. The first one is a semi-discrete cell-centered HTC finite volume scheme that employs collocated grids and is compatible with the total energy conservation law. However, these schemes generally do not respect the basic vector calculus identities at the discrete level. The second method is a staggered fully-discrete semi-implicit scheme. Thanks to a suitable vertex-based mesh staggering and via the use of compatible mimetic finite difference operators, exact discrete compatibility is achieved with the two classical vector calculus identities, namely $\nabla \cdot \nabla \times \A = 0$ and $\nabla \times \nabla \phi = 0$, for any generic vector and scalar field $ \phi $ and $  \A $, respectively. The second finite volume scheme proposed in this paper can be also proven to be total energy conservative at the fully-discrete level.
We applied both numerical schemes to three systems: nonlinear acoustic equations, nonlinear Maxwell equations in the absence of charges, and the nonlinear Maxwell-GLM equations. In all cases, we have proven total energy conservation at the discrete level. Furthermore, we have demonstrated that the divergence and curl type involutions are precisely satisfied at the discrete level when the second fully discrete method is employed.

\section*{Acknowledgments}

This research was funded by the European Union NextGenerationEU project PNRR Spoke 7 CN HPC.  Views and opinions expressed are however those of the author(s) only and do not necessarily reflect those of the European Union or the European Research
Council. Neither the European Union nor the granting authority can be held responsible for them. 

A.L. and M.D. are members of the INdAM GNCS group and also acknowledge the financial support received from  
the Italian Ministry of Education, University and Research (MIUR) in the frame of the PRIN 2022 project \textit{High order structure-preserving semi-implicit schemes for hyperbolic equations}. M.D. also received funding via the Fondazione Caritro under the project SOPHOS. 

\section*{Data Availability}

The data can be obtained from the authors on reasonable request. 

\section*{Conflict of interest} 
 The authors declare that they have no conflict of interest. 

\bibliographystyle{plain}
\bibliography{biblio}
\end{document}